\pdfoutput=1
\RequirePackage{ifpdf}
\ifpdf % We~are running pdfTeX in pdf mode
\documentclass[pdftex]{sigma}
\else
\documentclass{sigma}
\fi

\numberwithin{equation}{section}

\newtheorem{Theorem}{Theorem}[section]
\newtheorem{Corollary}[Theorem]{Corollary}
\newtheorem{Lemma}[Theorem]{Lemma}
\newtheorem{Proposition}[Theorem]{Proposition}
 { \theoremstyle{definition}
\newtheorem{Definition}[Theorem]{Definition}
\newtheorem{Example}[Theorem]{Example}
\newtheorem{Remark}[Theorem]{Remark} }

\begin{document}
\allowdisplaybreaks

\newcommand{\arXivNumber}{1811.05084}

\renewcommand{\PaperNumber}{121}

\FirstPageHeading

\ShortArticleName{Obstructions for Symplectic Lie Algebroids}

\ArticleName{Obstructions for Symplectic Lie Algebroids}

\Author{Ralph L.~KLAASSE}

\AuthorNameForHeading{R.L.~Klaasse}

\Address{D\'epartement de Mathematique, Universit\'e libre de Bruxelles,\\ CP 218 Boulevard du Triomphe, B-1050 Bruxelles, Belgium}
\Email{\href{mailto:r.l.klaasse@gmail.com}{r.l.klaasse@gmail.com}}
\URLaddress{\url{http://homepages.ulb.ac.be/~rklaasse/}}

\ArticleDates{Received April 06, 2020, in final form November 23, 2020; Published online November 27, 2020}

\Abstract{Several types of generically-nondegenerate Poisson structures can be effectively studied as symplectic structures on naturally associated Lie algebroids. Relevant examples of this phenomenon include log-, elliptic, $b^k$-, scattering and elliptic-log Poisson structures. In this paper we discuss topological obstructions to the existence of such Poisson structures, obtained through the characteristic classes of their associated symplectic Lie algebroids. In particular we obtain the full obstructions for surfaces to carry such Poisson structures.}

\Classification{53D17; 53D05}

\section{Introduction}\label{sec:introduction}

Generically-nondegenerate Poisson structures have recently seen an intense increase in interest. The main reason for this has been the ability to effectively study them using Lie algebroids. Namely, in several instances it is possible to, given a Poisson structure $\pi \in {\rm Poiss}(X)$, define a~Lie algebroid $\mathcal{A} \to X$ adhering to the same mild degeneracies as $\pi$, such that $\pi$ is in a precise sense dual to a symplectic structure in $\mathcal{A}$, i.e., a closed nondegenerate $\mathcal{A}$-two-form.

Symplectic Lie algebroids were first considered in \cite{NestTsygan01}, and have more recently been studied especially when the anchor map $\rho_\mathcal{A}\colon \mathcal{A} \to TX$ is generically an isomorphism. This class includes log- \cite{Cavalcanti17,GualtieriLi14,GuilleminMirandaPires14,MarcutOsornoTorres14,MarcutOsornoTorres14two}, elliptic \cite{CavalcantiGualtieri18, CavalcantiKlaasse18}, $b^k$- \cite{GuilleminMirandaWeitsman17, MirandaPlanas18two, MirandaPlanas18, Scott16} and scattering symplectic structures \cite{Lanius16}. Through the use of symplectic Lie algebroids, powerful symplectic techniques can be brought to bear to study the associated Poisson structures, leading to various results.

In this paper we are interested in obtaining obstructions to the existence of a symplectic structure on a Lie algebroid, and thus to their underlying Poisson structures. While we focus in this paper primarily on Lie algebroids and their symplectic structures, these should be thought of as tools to make statements about interesting classes of Poisson structures.

A symplectic manifold inherits a natural orientation, and is further almost-complex. Analogous statements hold for any symplectic Lie algebroid $\mathcal{A} \to X$ (see Proposition~\ref{prop:symplaconsequences}), or indeed any symplectic vector bundle (without integrability condition). The existence of an orientation and complex structure on $\mathcal{A}$ is determined by the underlying vector bundle, and is obstructed by its characteristic classes. Much of this paper consists of the computation of characteristic classes for several specific Lie algebroids $\mathcal{A}$. Before we can state our results (Theorems~\ref{thm:introorientable} and~\ref{thm:intrologsymp}) we must first recall how these Lie algebroids can be constructed.

\subsection{Lie algebroids from divisors}

One of the tools we use is the language of (real) divisors on smooth manifolds \cite{CavalcantiKlaasse18,Klaasse18}. A \emph{$($real$)$ divisor} on~$X$ is a pair $(U,\sigma)$ consisting of a real line bundle with a section $\sigma \in \Gamma(U)$ that has nowhere dense zero set $Z_\sigma = \sigma^{-1}(0)$. Through evaluation via the map $\sigma\colon \Gamma(U^*) \to C^\infty(X)$ we obtain a \emph{divisor ideal} $I_\sigma \subseteq C^\infty(X)$. A divisor ideal determines a divisor up to line bundle isomorphism and multiplication by a nonvanishing smooth function, allowing us to mostly work with divisor ideals.

In this paper we will use the following three examples of divisors (see \cite{CavalcantiGualtieri18,CavalcantiKlaasse18,Klaasse18}):
\begin{itemize}\itemsep=0pt
\item \emph{log divisors}, denoted as $(L,s)$, where $s$ has transverse zeroes. Here $Z := Z_s$ is a co\-di\-men\-sion-one hypersurface, and $I_Z := I_s$ is its vanishing ideal, locally we have $I_Z = \langle z \rangle$ with $Z = \{z = 0\}$, so that the log divisor (and $I_Z$) is determined by $Z$.
\item \emph{elliptic divisors}, denoted as $(R,q)$, where $q$ along its smooth codimension-two zero set $D := Z_q$ has definite Hessian ${\rm Hess}(q) \in \Gamma\big(D;{\rm Sym}^2 N^*D \otimes R\big)$. Its divisor ideal is denoted $I_{|D|} := I_q$, and is locally given by $I_{|D|} = \big\langle r^2 \big\rangle$ with $r$ a radial distance in~$ND$;
\item \emph{elliptic-log divisors} $(L,s) \otimes (R,q)$, obtained as the product of a log and elliptic divisor such that $D \subseteq Z$. Its divisor ideal is $I_W = I_Z \cdot I_{|D|}$, and locally $I_W = \big\langle r^3 \cos \theta \big\rangle$.
\end{itemize}
Note that elliptic and elliptic-log divisors are not determined by their underlying zero sets. Regardless, we will write $(X,Z)$, $(X,|D|)$ and $(X,W)$ for \emph{log}, \emph{elliptic}, and \emph{elliptic-log pairs} for manifold pairs equipped with the above three divisor types respectively, and denote by $L_Z$ the line bundle of the log divisor associated to the pair $(X,Z)$.

An immediate consequence of the definition is the following, which we will use later.

\begin{Lemma}\label{lem:logdivisor} Let $(X,Z)$ be a log pair. Then we have $w_1(L_Z) = {\rm PD}_{\mathbb{Z}_2}[Z] \in H^1(X;\mathbb{Z}_2)$.
\end{Lemma}

Here $w_1$ is the first Stiefel--Whitney class, and ${\rm PD}_{\mathbb{Z}_2}$ is the Poincar\'e dual with $\mathbb{Z}_2$-coefficients. This follows because the section $s \in \Gamma(L_Z)$ can be used to determine the Euler class of $L_Z$, which equals the top Stiefel--Whitney class (in this case $w_1$) after reduction modulo two.

\begin{Remark}\label{rem:separating} Note that ${\rm PD}_{\mathbb{Z}_2}[Z] = 0$ if and only if $Z$ separates $X$. That is, if and only if $Z$ has a \emph{global defining function} $f\colon X \to \mathbb{R}$, with $0$ a regular value of $f$ and $Z = f^{-1}(0)$.
\end{Remark}

The next step is to recall that a \emph{Lie algebroid} is a vector bundle $\mathcal{A} \to X$ equipped with an \emph{anchor} map $\rho_{\mathcal{A}}\colon \mathcal{A} \to TX$ and a Lie bracket $[\cdot,\cdot]_\mathcal{A}$ on $\Gamma(\mathcal{A})$ which satisfies $[v, f w]_\mathcal{A} = f [v,w]_\mathcal{A} + \mathcal{L}_{\rho_{\mathcal{A}}(v)} f \cdot w$ for all $v, w \in \Gamma(\mathcal{A})$ and $f \in C^\infty(X)$. Divisor ideals are an effective tool to construct Lie algebroids generically isomorphic to $TX$, as we now explain.

Let $I \subseteq C^\infty(X)$ be a divisor ideal and $\Gamma(TX)_I = \{V \in \Gamma(TX)\colon \mathcal{L}_V I \subseteq I\} \subseteq \Gamma(TX)$ be the involutive submodule of vector fields preserving $I$. When $\Gamma(TX)_I$ is projective it specifies uniquely a Lie algebroid $\mathcal{A}_I \to X$ such that $\Gamma(\mathcal{A}_I) \cong \Gamma(TX)_I$ by the Serre--Swan theorem.

\begin{Definition}[\cite{Klaasse18}] Let $I \subseteq C^\infty(X)$ be a divisor ideal for which $\Gamma(TX)_I$ is projective. Then the Lie algebroid $\mathcal{A}_I \to X$ with $\Gamma(\mathcal{A}_I) \cong \Gamma(TX)_I$ is called the \emph{ideal Lie algebroid} of $I$.
\end{Definition}

Examples of this construction include:
\begin{itemize}\itemsep=0pt
\item The \emph{log-tangent bundle} $\mathcal{A}_Z = TX(-\log Z)$ associated to $I_Z$ \cite{Melrose93};
\item The \emph{elliptic tangent bundle} $\mathcal{A}_{|D|} = TX(-\log |D|)$ associated to $I_{|D|}$ \cite{CavalcantiGualtieri18};
\item The \emph{elliptic-log tangent bundle} $\mathcal{A}_W = TX(-\log W)$ associated to $I_W$ \cite{Klaasse18}.	
\end{itemize}
Note that the latter has natural morphisms onto $\mathcal{A}_D$ and $\mathcal{A}_Z$ via the section module inclusion. These Lie algebroids all have the property that their anchor $\rho_\mathcal{A}\colon \mathcal{A} \to TX$ is an isomorphism on a dense open set, which is the complement of their \emph{degeneracy locus}. For these Lie algebroids the anchor map gives a divisor ${\rm div}(\mathcal{A}) = (\det(TX) \otimes \det(\mathcal{A}^*), \det(\rho_\mathcal{A}))$ with divisor ideal $I_\mathcal{A}$.

One can also obtain Lie algebroids by modifying a given Lie algebroid using a Lie subalgebroid supported on a hypersurface $Z$. This process is called \emph{(lower) elementary modification} \cite{GualtieriLi14,Li13} or \emph{rescaling} \cite{Lanius16, Melrose93}. This can be extended to divisor ideals $I \subseteq C^\infty(X)$ supported on smooth submanifolds other than log ideals $I_Z$ (see \cite{Klaasse18}), but we will not use this here.

Before we can provide the definition, recall that a \emph{Lie subalgebroid} of $\mathcal{A} \to X$ is a Lie algebroid $\mathcal{B} \to Z$ with $Z \subseteq X$ carrying an injective Lie algebroid morphism covering an injective immersion: $(\varphi,i)\colon (\mathcal{B},Z) \hookrightarrow (\mathcal{A},X)$. Here a Lie algebroid morphism over varying base is most succinctly defined as a vector bundle morphism $\varphi\colon \mathcal{B} \to \mathcal{A}$ intertwining their differentials, i.e., such that $d_\mathcal{B} \circ \varphi^* = \varphi^* \circ d_\mathcal{A}$.

\begin{Definition} Let $(X,Z)$ be a log pair and $(\mathcal{B},Z) \subseteq (\mathcal{A},X)$ a Lie subalgebroid. The \emph{lower elementary modification} or \emph{$(\mathcal{B},Z)$-rescaling} of $\mathcal{A}$ along $\mathcal{B}$ is the Lie algebroid $[\mathcal{A}\text{:}\mathcal{B}]$ defined by
	\begin{displaymath}
	\Gamma([\mathcal{A}\text{:}\mathcal{B}]) \cong \{v \in \Gamma(\mathcal{A}) \colon v|_Z \in \Gamma(\mathcal{B})\}.
	\end{displaymath}
\end{Definition}

Elementary modification can also be performed purely at the level of vector bundles.

\begin{Remark}\label{rem:zerorescaling} Given a Lie algebroid $\mathcal{A} \to X$ and a hypersurface $Z \subseteq X$, one can always perform $(0,Z)$-rescaling. The resulting Lie algebroid $[\mathcal{A}\text{:}0]$ is isomorphic to the tensor product $\mathcal{A} \otimes L_Z$ as a vector bundle.
\end{Remark}

\begin{Example} Let $(X,Z)$ be a log pair. The following are examples of modifications \cite{Melrose93,Melrose95}:
	\begin{itemize}\itemsep=0pt
	\item the log-tangent bundle $\mathcal{A}_Z = [TX\text{:}TZ]$, locally given by $\Gamma(\mathcal{A}_Z) = \langle z \partial_z, \partial_{x_i} \rangle$;
	\item the \emph{zero tangent bundle} $\mathcal{B}_Z = [TX\text{:}0]$, locally given by $\Gamma(\mathcal{B}_Z) = \langle z \partial_z, z \partial_{x_i} \rangle$;
	\item the \emph{scattering tangent bundle} $\mathcal{C}_Z = [\mathcal{A}_Z\text{:}0]$, locally given by $\Gamma(\mathcal{C}_Z) = \langle z^2 \partial_z, z \partial_{x_i} \rangle$.
	\end{itemize}
\end{Example}

Given a log pair $(X,Z)$ and $k \geq 1$ we can define a Lie algebroid $\mathcal{A}_Z^k \to X$ as follows. Using the inclusion $\iota_Z\colon Z \hookrightarrow X$ and the vanishing ideal sheaf $\mathcal{I}_Z$ for $Z$, consider the sheaf $\mathcal{J}_Z^k := \iota_Z^{-1}\big(C^\infty_X / \mathcal{I}_Z^{k+1}\big)$. Denote its space of global sections by~$J_Z^k$, and fix $j \in J_Z^{k-1}$, which is called a $(k-1)$-jet. Given a smooth function $f$ defined in a neighbourhood of~$Z$, we write $f \in j$ if $f$ represents $j$. Assume that $j$ is represented by local defining functions for~$Z$, making it a \emph{defining $(k-1)$-jet} (we suppress this from our notation). With this, we define $\mathcal{A}_Z^k$ by setting
\begin{displaymath}
\Gamma\big(\mathcal{A}_Z^k\big) \cong \big\{V \in \Gamma(TX)\colon \mathcal{L}_V f \in I_Z^k \text{ for all } f \in j\big\}.
\end{displaymath}
This is the \emph{$b^k$-tangent bundle}~\cite{Scott16}, and is locally given by $\Gamma\big(\mathcal{A}_Z^k\big) = \big\langle z^k \partial_z, \partial_{x_i} \big\rangle$ for a local $z \in j$. When $k = 1$ the jet data is vacuous, so that $\mathcal{A}_Z^1 = \mathcal{A}_Z$, the log-tangent bundle.

\subsection{Poisson structures on Lie algebroids} Poisson structures are readily linked to divisors and the Lie algebroids built from them (see~\cite{Klaasse18}). Let $\pi \in {\rm Poiss}\big(X^{2n}\big)$ be a Poisson structure, and consider its \emph{Pfaffian}, $\wedge^n \pi \in \Gamma(\det(TX))$. If $\pi$ is generically nondegenerate this defines a divisor $(\det(TX), \wedge^n \pi)$ and hence a divisor ideal $I_\pi$. We say $\pi$ is of \emph{$I$-divisor-type} if $I_\pi = I$.

Poisson structures on a Lie algebroid $\mathcal{A} \to X$ are defined as those sections $\pi_\mathcal{A} \in \Gamma\big({\wedge}^2 \mathcal{A}\big)$ such that $[\pi_\mathcal{A},\pi_\mathcal{A}]_\mathcal{A} = 0$. These specify underlying Poisson structures $\pi = \rho_\mathcal{A}(\pi_\mathcal{A}) \in {\rm Poiss}(X)$. In~\cite{Klaasse18} we showed that if $\pi \in {\rm Poiss}(X)$ is of $I$-divisor-type, and $I$ is such that its ideal Lie algebroid $\mathcal{A}_I$ exists, then $\pi$ admits an \emph{$\mathcal{A}_I$-lift}: there exists a (unique) $\mathcal{A}_I$-Poisson structure $\pi_{\mathcal{A}_I}$ such that $\pi = \rho_{\mathcal{A}_I}(\pi_{\mathcal{A}_I})$. We say that a divisor ideal $I$ is \emph{standard} if its ideal Lie algebroid satisfies $I_{\mathcal{A}_I} = I$. As noted in \cite{Klaasse18}, log, elliptic, and elliptic-log divisor ideals are standard. If the divisor ideal $I$ is standard, then the lifted Poisson structure $\pi_{\mathcal{A}_I}$ is nondegenerate.

A Lie algebroid $\mathcal{A} \to X$ of even rank is \emph{symplectic} if it carries a nondegenerate closed $\mathcal{A}$-two-form $\omega_\mathcal{A}$ (after~\cite{NestTsygan01}). Such an \emph{$\mathcal{A}$-symplectic structure} corresponds to a nondegenerate $\mathcal{A}$-Poisson structure~$\pi_{\mathcal{A}}$ via the relation $\pi_\mathcal{A}^\sharp = \big(\omega_\mathcal{A}^\flat\big)^{-1}$. Due to this, if we wish to study Poisson structures of $I$-divisor-type, we can study $\mathcal{A}_I$-symplectic geometry instead.

Summarizing the above definitions, if we have a log pair $(X,Z)$, an elliptic pair $(X,|D|)$ or an elliptic-log pair $(X,W)$, we can define and study the following classes of Poisson structures and symplectic Lie algebroids:
\begin{itemize}\itemsep=0pt
\item \emph{log-Poisson structures}, where $\pi$ is of $I_Z$-divisor-type. These are also \emph{log-symplectic structures}, associated to the Lie algebroid $\mathcal{A}_Z$ \cite{GuilleminMirandaPires14, MarcutOsornoTorres14two}, also \cite{Cavalcanti17,GualtieriLi14, MarcutOsornoTorres14} and others;
\item \emph{elliptic Poisson structures}, where $\pi$ is of $I_{|D|}$-divisor-type. These are also \emph{elliptic symplectic structures}, associated to the Lie algebroid $\mathcal{A}_{|D|}$ \cite{CavalcantiGualtieri18}, also \cite{CavalcantiKlaasse18};
\item \emph{elliptic-log Poisson structures}, where $\pi$ is of $I_{W}$-divisor-type. These are related to \emph{elliptic-log symplectic structures}, associated to the Lie algebroid $\mathcal{A}_W$ \cite{Klaasse18};
\end{itemize}
More directly defined through their symplectic Lie algebroids, we have:
\begin{itemize}\itemsep=0pt
\item \emph{zero symplectic structures}, associated to $\mathcal{B}_Z$ (cf.\ \cite{Lanius16}, and Remark~\ref{rem:zeronosymp});
\item \emph{scattering symplectic structures}, associated to $\mathcal{C}_Z$ \cite{Lanius16};
\item \emph{$b^k$-symplectic structures}, associated to $\mathcal{A}_Z^k$ \cite{Scott16}, also \cite{GuilleminMirandaWeitsman17}.
\end{itemize}
Each of these has an underlying Poisson structure, which can often be characterized intrinsically. It is not our intent to describe the geometry of these Poisson structures in great detail here. While the remainder of this paper uses Lie algebroids and Lie algebroid objects, these are viewed as tools to make statements about generically-nondegenerate Poisson structures.

With this in mind we can state our results.

\subsection{Results}

Our first result is the following (Theorem~\ref{thm:orientable}), regarding orientability. Recall that a mani\-fold~$X$ is orientable if and only if the first Stiefel--Whitney class of its tangent bundle vanishes, i.e., $w_1(TX) = 0 \in H^1(X;\mathbb{Z}_2)$. The following are analogues of this.

\begin{Theorem}\label{thm:introorientable} Let $\mathcal{A} \to X^n$ be a symplectic Lie algebroid. Then in $H^1(X;\mathbb{Z}_2)$ we have:
\begin{itemize}\itemsep=0pt
	\item $w_1(TX) + k {\rm PD}_{\mathbb{Z}_2}[Z] = 0$ if $\mathcal{A} = \mathcal{A}_Z^k$, the $b^k$-tangent bundle;
	\item $w_1(TX) = 0$ if $\mathcal{A} = \mathcal{B}_Z$, the zero tangent bundle;
	\item $w_1(TX) + {\rm PD}_{\mathbb{Z}_2}[Z] = 0$ if $\mathcal{A} = \mathcal{C}_Z$, the scattering tangent bundle;
	\item $w_1(TX) = 0$ if $\mathcal{A} = \mathcal{A}_{|D|}$, the elliptic tangent bundle;
	\item $w_1(TX) + {\rm PD}_{\mathbb{Z}_2}[Z] = 0$ if $\mathcal{A} = \mathcal{A}_W$, the elliptic-log tangent bundle.
	\end{itemize}
Here $w_1$ is the first Stiefel--Whitney class, and ${\rm PD}_{\mathbb{Z}_2}$ is the Poincar\'e dual with $\mathbb{Z}_2$-coefficients.
\end{Theorem}

This result provides the full obstruction for a surface to be $\mathcal{A}$-symplectic. This latter statement is because the integrability condition (closedness) is immediate, so that only a nondegenerate $\mathcal{A}$-two-form is required, which exists if and only if $\mathcal{A}$ satisfies $w_1(\mathcal{A}) = 0$, i.e., if and only if its first Stiefel--Whitney class vanishes in $H^1(X;\mathbb{Z}_2)$.

\begin{Remark} If we combine Theorem~\ref{thm:introorientable} with Remark~\ref{rem:separating}, we see that the singular locus $Z$ of a log-symplectic or scattering-symplectic manifold $X$ admits a global defining function if and only if $X$ is orientable, and similarly for $b^k$-symplectic structures when $k$ is odd.
\end{Remark}

Our second result draws consequences from the required complex structure on the Lie algebroids in dimension four (see Theorems~\ref{thm:azsympstr} and~\ref{thm:scatteringsymp}). It is the analogue of Noether's formula \cite[Theorem~1.4.13]{GompfStipsicz99} that exists for regular symplectic four-manifolds.

\begin{Theorem}\label{thm:intrologsymp} Let $\big(X^4,Z\big)$ be a compact oriented $b^k$-symplectic or scattering-symplectic four-manifold, with $k$ odd. Then $b_2^+(X) + b_1(X) + f(X,Z)$ is odd, where $2 f(X,Z) = e(\mathcal{A}_Z) - e(TX)$.
\end{Theorem}

In the above, $b_1(X)$ is the first Betti number of $X$, and $b_2^+(X)$ is the dimension of a maximal positive definite subspace on $H^2(X;\mathbb{R})$ with respect to the natural quadratic form present on $H^2(X;\mathbb{R})$ in dimension four. Finally, for an oriented vector bundle $E^n \to X^n$, its Euler class is denoted by $e(E) \in H^n(X;\mathbb{Z})$.

\section{Computing characteristic classes}\label{sec:homotopical}

We start towards obtaining homotopical obstructions to the existence of~$\mathcal{A}$-symplectic structures on a given closed manifold~$X$. More precisely, we focus on the following simple facts.

\begin{Proposition}\label{prop:symplaconsequences} Let $\mathcal{A} \to X$ be a symplectic Lie algebroid. Then:
	\begin{itemize}\itemsep=0pt
	\item $\mathcal{A}$ must be orientable, i.e., it must satisfy $w_1(\mathcal{A}) = 0 \in H^1(X;\mathbb{Z}_2)$;
	\item $\mathcal{A}$ must be complex, i.e., there must exist a $J_\mathcal{A} \in {\rm End}(\mathcal{A})$ with $J_\mathcal{A}^2 = - {\rm id}$.
	\end{itemize}
\end{Proposition}

These properties both follow from the linear algebra of having a nondegenerate $\mathcal{A}$-two-form. Indeed, they hold for any symplectic vector bundle (e.g.,~\cite{McDuffSalamon17}), as they do not use integrability.

\begin{proof} Let $\omega_\mathcal{A}$ be an $\mathcal{A}$-symplectic structure. Then ${\rm rank}(\mathcal{A}) = 2m$ is necessarily even, and $\omega_{A}^m \in \Gamma(\det(\mathcal{A}^*))$ is nonvanishing. Thus $\det(\mathcal{A}^*)$ is trivial, and $w_1(\mathcal{A}) = w_1(\det(\mathcal{A}^*)) = 0$. To see $\mathcal{A}$ must admit a complex structure, follow the standard proof for $\mathcal{A} = TX$ (e.g.,~\cite{McDuffSalamon17}).
\end{proof}

Note that when both $\mathcal{A}$ and $X$ are four-dimensional, a classical result by Wu~\cite{Wu52} (see also~\cite{HirzebruchHopf58}) can be used, characterizing when an oriented vector bundle admits a complex structure.

\begin{Theorem}[\cite{Wu52}]\label{thm:wuacs} Let $E^4 \to X^4$ be an oriented Euclidean rank-four vector bundle over a~compact oriented four-manifold. Then $E$ admits a complex structure if and only if there exists a class $c \in H^2(X;\mathbb{Z})$ such that $c \!\!\mod 2 \equiv w_2(E) \in H^2(X;\mathbb{Z}_2)$ and $c^2 = p_1(E) + 2 e(E)$.
\end{Theorem}

To make effective use of these observations, it is clear that we must determine the relevant characteristic classes of the bundle $\mathcal{A} \to X$. We do this via stable bundle isomorphisms.

\subsection{Stable bundle isomorphisms}\label{sec:stablebundleiso}

For the $b^k$-tangent bundles we can establish a stable bundle isomorphism, relating $\mathcal{A}_Z^k$ to $TX$. Denote by ${\underline{\mathbb{R}}} \to X$ the trivial real line bundle.

\begin{Proposition}\label{prop:bkbundleiso} Let $\big(X^n,Z\big)$ be a log pair with a defining $(k-1)$-jet $j$ at $Z$. Then we have $\mathcal{A}_Z^k \oplus L_Z \cong TX \oplus \underline{\mathbb{R}}$ if $k$ is odd. Moreover, for any $k \geq 3$ we have $\mathcal{A}_Z^k \oplus \underline{\mathbb{R}} \cong \mathcal{A}_Z^{k-2} \oplus \underline{\mathbb{R}}$.
\end{Proposition}

We emphasize that these are vector bundle isomorphisms, and not of Lie algebroids.

\begin{Remark}\label{rem:slightinaccuracy} When $k = 1$, Proposition~\ref{prop:bkbundleiso} reduces to the statement that
\begin{displaymath}
	\mathcal{A}_Z \oplus L_Z \cong TX \oplus \underline{\mathbb{R}}.
\end{displaymath}
This was noted in \cite{CannasDaSilva10, CannasDaSilvaGuilleminWoodward00} when $X$ is orientable without proof, and they state that in general one has $\mathcal{A}_Z \oplus \underline{\mathbb{R}} \cong TX \oplus L_Z$. When $Z$ is separating, so that $L_Z$ is trivial (see Lemma~\ref{lem:logdivisor}), these two statements are equivalent; for example, when $(X,Z)$ is log-symplectic and $X$ is orientable.
\end{Remark}

\begin{proof} From the $(k-1)$-jet $j$ at $Z$, we can choose a compatible log divisor $(L_Z,s)$ such that in any local trivialization, the section $s$ represents $j$. We will first prove the first assertion, hence assume that $k$ is odd. Over $X \backslash Z$, the line bundle $L_Z$ is trivial, and moreover $\mathcal{A}_Z^k \cong TX$ via the anchor map $\rho$, so that there we have a clear choice of isomorphism $\varphi\colon \mathcal{A}_Z^k \oplus L_Z \to TX \oplus \underline{\mathbb{R}}$. In matrix form $\varphi$ becomes as follows, using the evaluation $\langle s , - \rangle\colon \Gamma(L_Z) \to \Gamma(\underline{\mathbb{R}})$ induced by $s$:
\begin{displaymath}
 	\varphi = \begin{pmatrix} \rho & 0 \\ 0 & \langle s , - \rangle \end{pmatrix}.
\end{displaymath}
Near $Z$ we define a bundle isomorphism extending the one above as follows.

Choose a tubular neighbourhood embedding of $Z$. Then, within its image (which we suppress from our notation), choose a trivializing open cover $\{U_\alpha\}_\alpha$ around $Z$ on which $L_Z$ is trivial with nonvanishing sections $\tau_\alpha \in \Gamma(L_Z|_{U_\alpha})$. We then shrink the open cover so that the associated opens $V_\alpha := U_\alpha \cap Z$ on $Z$ are such that $TZ|_{V_\alpha}$ is also trivial.

The given section $s \in \Gamma(L_Z)$ vanishes transversally, which by definition means that its normal derivative provides an isomorphism $d^\nu s\colon NZ \stackrel{\cong}{\to} L_Z|_Z$. Because of this, if we use the resulting trivialization for $NZ$ on $V_\alpha$, we have now that $U_\alpha \cong V_\alpha \times (-1,1)$, say, with the coordinate $z_\alpha$ of the second factor also being such that $s_\alpha = z_\alpha \tau_\alpha$. We define the transition functions for $L_Z$ by the relation $\tau_\alpha = g_\alpha^\beta \tau_\beta$ on $U_\alpha \cap U_\beta$. Due to these choices, then, we have that
\begin{displaymath}
	\Gamma\big(U_\alpha; \mathcal{A}_Z^k\big) = \big\langle z_\alpha^k \partial_{z_\alpha}, \{v_{\alpha,i}\}_{i=2}^n\big\rangle, \qquad \text{and} \qquad \Gamma(U_\alpha;TX) = \big\langle \partial_{z_\alpha}, \{v_{\alpha,i}\}_{i=2}^n \big\rangle,
\end{displaymath}
with $z_\alpha \in j$ in the given $(k-1)$-jet. Generic sections of $\mathcal{A}_Z^k \oplus L_Z$ and $TX \oplus \underline{\mathbb{R}}$ on $U_\alpha$ can then expressed respectively as the tuples
\begin{displaymath}
	\bigg(\sum_{2 \leq i \leq n} \lambda_i \cdot v_{\alpha,i} + \lambda_1 \cdot z_\alpha^k \partial_{z_\alpha}, \lambda_{n+1} \cdot \tau_\alpha\bigg) \qquad \text{and} \qquad \bigg(\sum_{2\leq i \leq n} \mu_i \cdot v_{\alpha,i} + \mu_1 \cdot \partial_{z_\alpha}, \mu_{n+1} \cdot \mathbf{1}\bigg),
\end{displaymath}
where $\lambda_i, \mu_i \in C^\infty(U_\alpha)$. Choose a bundle metric on $NZ$ and a radial bump function $\psi$ whose support is contained in a disk bundle around $Z$ of the tubular neighbourhood embedding, and which equals $1$ on $Z$. Denote by $A_\alpha\colon \Gamma(U_\alpha;L_Z) \to \Gamma(U_\alpha;TX)$ the map sending $\tau_\alpha$ to $\partial_{z_\alpha}$, and similarly let $B_\alpha\colon \Gamma\big(U_\alpha;\mathcal{A}_Z^k\big) \to \Gamma(U_\alpha; \underline{\mathbb{R}})$ be the map sending $z_\alpha^k \partial_{z_\alpha}$ to $\mathbf{1}$.

In the given tubular neighbourhood we now define the map $\varphi\colon \mathcal{A}_Z^k \oplus L_Z \to TX \oplus \underline{\mathbb{R}}$ on sections to be in matrix form given by
\begin{displaymath}
	\varphi = \begin{pmatrix} \rho & \psi A_\alpha \\ -\psi B_\alpha & \langle s , - \rangle \end{pmatrix} = \begin{pmatrix} I_{n-1} & 0 & 0 \\ 0 & z_\alpha^k & \psi \\ 0 & -\psi & z_\alpha \end{pmatrix},
\end{displaymath}
where $I_{n-1}$ is the $(n-1) \times (n-1)$ identity matrix. The maps $B_\alpha$ are well-defined because the local sections $z_\alpha^k \partial_{z_\alpha}$ together define a nowhere-vanishing section in the kernel of $\rho|_Z\colon \mathcal{A}_Z^k|_Z \to TM|_Z$ as described in \cite[Proposition 2.3]{Scott16} (for $k = 1$ this can also be found in \cite{GuilleminMirandaPires14} and follows directly). For well-definedness of the maps $A_\alpha$, note that on intersections $U_\alpha \cap U_\beta$ we have that
\begin{displaymath}
	s_\alpha = z_\alpha \tau_\alpha = z_\alpha g_\alpha^\beta \tau_\beta = z_\beta \tau_\beta = s_\beta,
\end{displaymath}
so that $z_\alpha g_\alpha^\beta = z_\beta$, from which $\partial_{z_\alpha} = g_\alpha^\beta \partial_{z_\beta}$ follows, and hence well definedness is clear:
\begin{displaymath}
	\partial_{z_\alpha} = A_\alpha(\tau_\alpha) = A_\beta\big(g_\alpha^\beta \tau_\beta\big)= g_\alpha^\beta \partial_{z_\beta}.
\end{displaymath}
To see that $\varphi$ is an isomorphism, we merely note that in matrix form the map $\varphi$ has determinant equal to $z_\alpha^{k+1} + \psi^2$ on $U_\alpha$, which is always nonvanishing because $k$ is odd, showing invertibility. It is clear that the thus defined bundle map $\varphi$ extends the one on $X \backslash Z$ described before.

Next, let $k \geq 3$ be given and consider the Lie algebroid morphism $\rho_k^{k-2} \colon \mathcal{A}_Z^k \to \mathcal{A}_Z^{k-2}$ induced from the natural map of viewing a $(k-1)$-jet as a $(k-3)$-jet. We use the same trivializations and bump function $\psi$ as above, and define a map $\varphi'\colon \mathcal{A}_Z^k \oplus \underline{\mathbb{R}} \to \mathcal{A}_Z^{k-2} \oplus \underline{\mathbb{R}}$ in matrix form by
\begin{displaymath}
	\varphi' = \begin{pmatrix} \rho_k^{k-2} & \psi C_\alpha \\ -\psi D_\alpha & {\rm id} \end{pmatrix} = \begin{pmatrix} I_{n-1} & 0 & 0 \\ 0 & z_\alpha^2 & \psi \\ 0 & -\psi & 1 \end{pmatrix},
\end{displaymath}
with $C_\alpha\colon \underline{\mathbb{R}} \to \mathcal{A}_Z^{k-2}$ the map $\mathbf{1} \mapsto z_\alpha^{k-2} \partial_{z_\alpha}$ and $D_\alpha\colon \mathcal{A}_Z^k \to \underline{\mathbb{R}}$ the map $z_\alpha^k \partial_{z_\alpha} \to \mathbf{1}$. These maps are both well-defined as these four sections are global. The determinant of $\varphi'$ is seen to equal $z_\alpha^2 + \psi^2$ in these local coordinates, which is always non-vanishing, showing invertibility.
\end{proof}

\begin{Remark} One can not readily change this proof to instead show that $\mathcal{A}_Z \oplus \underline{\mathbb{R}} \cong TX \oplus L_Z$, as per Remark~\ref{rem:slightinaccuracy}, as in general there is no nonzero map from $\mathcal{A}_Z$ to $L_Z$. The crucial ingredient in the proof is the existence of the map from $L_Z$ to $TX$. The map $\varphi$ used in the proof also exists for $k$ even, but then is not an isomorphism.
\end{Remark}

A similar result holds more generally for $(\mathcal{B},Z)$-rescalings of $\mathcal{A}$ of higher corank, for which the quotient vector bundle $\mathcal{A}|_Z/\mathcal{B}$ is a sum of line bundles. This is proven by an induction-like adaptation of the same method, because we obtain a flag of subbundles $\mathcal{B} \subseteq \mathcal{B}_{k-1} \subseteq \dots \subseteq \mathcal{B}_1 \subseteq \mathcal{A}|_Z$ with ${\rm corank}(\mathcal{B}_i) = i$. This gives a bundle isomorphism similar to Proposition~\ref{prop:logbundleiso}.

\begin{Proposition}\label{prop:logbundleiso} Let $\mathcal{A} \to X$ be a vector bundle and $Z \subseteq X$ a hypersurface. Consider a $(\mathcal{B},Z)$-rescaling $[\mathcal{A}\text{:}\mathcal{B}]$ of $\mathcal{A}$ with ${\rm corank}(\mathcal{B}) = k$, and assume that the quotient bundle $\mathcal{A}|_Z/\mathcal{B} \to Z$ is a~sum of line bundles. Then using $k$ copies of $L_Z$ and $\underline{\mathbb{R}}$ we have
	\begin{displaymath}
	[\mathcal{A}\text{:}\mathcal{B}] \oplus L_Z \oplus \dots \oplus L_Z \cong \mathcal{A} \oplus \underline{\mathbb{R}} \oplus \dots \oplus \underline{\mathbb{R}}.
	\end{displaymath}
\end{Proposition}

We stress that the number of copies depends on the corank of the vector subbundle.

\subsection{Computing characteristic classes}\label{sec:compcharclasses}

In this section we compute relevant characteristic classes of the Lie algebroids we have introduced. We will mainly be interested in the first and second Stiefel--Whitney classes $w_1, w_2 \in H^{i}(X;\mathbb{Z}_2)$, and in the first Pontryagin class $p_1 \in H^4(X;\mathbb{Z})$. We recall several properties of these characteristic classes (see, e.g.,~\cite{MilnorStasheff74}).
\begin{Proposition}\label{prop:charclasses} Let $E^m, F^n \to X$ be real vector bundles. Denote the full Stiefel--Whitney and Pontryagin classes by $w\colon {\rm Vect}(X) \to H^\bullet(X;\mathbb{Z}_2)$ and $p\colon {\rm Vect}(X) \to H^\bullet(X;\mathbb{Z})$. Then:
	\begin{itemize}\itemsep=0pt
	\item[$i)$] $w(E \oplus F) = w(E) \cup w(F)$, and $w_1(E \otimes F) = n w_1(E) + m w_1(F)$;
	\item[$ii)$] $w_2(E \otimes F) = w_2(E) + w_1(F) \cup w_1(E)$ if $m = 4$ and $n = 1$.
	\item[$iii)$] $2 p(E \oplus F) = 2 p(E) \cup p(F)$, and $p(E \otimes F) = p(E)$ if $n = 1$;
	\end{itemize}
\end{Proposition}

We now determine the relevant characteristic classes for the Lie algebroids $\mathcal{A}_Z^k$, $\mathcal{B}_Z$, and $\mathcal{C}_Z$.

\begin{Proposition}\label{prop:logcharclass} Let $\big(X^n,Z\big)$ be a log pair with Lie algebroids $\mathcal{A}_Z^k$, $\mathcal{B}_Z$, and $\mathcal{C}_Z$, with $k$ odd. Then:
\begin{enumerate}\itemsep=0pt
\item[] for $\mathcal{A}^k_Z$:
\begin{itemize}\itemsep=0pt
\item $w_1\big(\mathcal{A}^k_Z\big) = w_1(TX) + w_1(L_Z)$,
\item $w_2\big(\mathcal{A}^k_Z\big) = w_2(TX) + w_1(L_Z) \cup w_1(TX)$,
\item $p_1\big(\mathcal{A}^k_Z\big) = p_1(TX)$ if $X$ is orientable and four-dimensional;
\end{itemize}
\item[]	for $\mathcal{B}_Z$:
\begin{itemize}\itemsep=0pt
\item $w_1(\mathcal{B}_Z) = w_1(TX) + n w_1(L_Z)$,
\item $w_2(\mathcal{B}_Z) = w_2(TX) + w_1(L_Z) \cup w_1(TX)$ if $X$ is four-dimensional;
\end{itemize}
\item[]	for $\mathcal{C}_Z$:
\begin{itemize}\itemsep=0pt
\item $w_1(\mathcal{C}_Z) = w_1(TX) + (n+1) w_1(L_Z)$,
\item $w_2(\mathcal{C}_Z) = w_2(TX)$ if $X$ is four-dimensional,
\item $p_1(\mathcal{C}_Z) = p_1(TX)$ if $X$ is orientable and four-dimensional.
\end{itemize}
\end{enumerate}
\end{Proposition}

\begin{Remark} The fact that $w(\mathcal{A}_Z) = w(TX)(1 + {\rm PD}_{\mathbb{Z}_2}[Z])$ can be found in \cite{CannasDaSilvaGuilleminWoodward00}, as this would also be what follows from the stable isomorphism relation noted by them, see Remark~\ref{rem:slightinaccuracy}.
\end{Remark}

\begin{proof}
{\it For $\mathcal{A}_Z^k$:} By Proposition~\ref{prop:bkbundleiso} we have $\mathcal{A}_Z^k \oplus L_Z \cong TX \oplus \underline{\mathbb{R}}$, hence due to Proposition~\ref{prop:charclasses}(i) we get $w\big(\mathcal{A}_Z^k\big)\cup (1 + w_1(L_Z)) = w(TX)$. In degree one this gives $w_1\big(\mathcal{A}_Z^k\big) = w_1(TX) + w_1(L_Z)$ as desired. In degree two it follows that $w_2\big(\mathcal{A}_Z^k\big) = w_2(TX) + w_1(L_Z) \cup w_1(TX)$. We see that $2 p_1\big(\mathcal{A}_Z^k\big) = 2 p_1(TX)$, as $p \equiv 1$ for line bundles. If $X$ is orientable and four-dimensional we know that $H^4(X;\mathbb{Z}) \cong \mathbb{Z}$, which in particular has no two-torsion, so that $p_1\big(\mathcal{A}_Z^k\big) = p_1(TX)$.

{\it For $\mathcal{B}_Z$:} This follows from Proposition~\ref{prop:charclasses} after using Remark~\ref{rem:zerorescaling} that $\mathcal{B}_Z \cong TX \otimes L_Z$.

{\it For $\mathcal{C}_Z$:} By Remark~\ref{rem:zerorescaling} we have $\mathcal{C}_Z \cong \mathcal{A}_Z \otimes L_Z$, so that from Proposition~\ref{prop:logbundleiso} for $\mathcal{A} = \mathcal{A}_Z$ we obtain $\mathcal{C}_Z \oplus L_Z^2 \cong \mathcal{B}_Z \oplus L_Z$. As $L_Z^2$ is canonically trivial, using Proposition~\ref{prop:charclasses}(i) this gives $w(\mathcal{C}_Z) = (1 + w_1(L_Z)) \cup w(\mathcal{B}_Z)$. In degree one this results in (using the case of $\mathcal{B}_Z$ above):
\begin{displaymath}
	w_1(\mathcal{C}_Z) = w_1(\mathcal{B}_Z) + w_1(L_Z) = w_1(TX) + (n+1)w_1(L_Z).
\end{displaymath}
In degree two we see similarly that if $X$ is four-dimensional that
\begin{align*}
w_2(\mathcal{C}_Z) &= w_2(\mathcal{B}_Z) + w_1(L_Z) \cup w_1(\mathcal{B}_Z) \\
&=w_2(TX) + w_1(L_Z) \cup w_1(TX) + w_1(L_Z) \cup (w_1(TX) + 4 w_1(L_Z)) = w_2(TX).
\end{align*}
Assuming also orientability of $X$, Proposition~\ref{prop:charclasses} and the case of $\mathcal{A}_Z$ determine $p_1(\mathcal{C}_Z)$.
\end{proof}

We can compute these characteristic classes somewhat more generally for rescalings.

\begin{Proposition} Let $\mathcal{A} \to X$ be a vector bundle and $Z \subseteq X$ a hypersurface. Let $[\mathcal{A}\text{:}\mathcal{B}] \to X$ be a corank-$k$ $(\mathcal{B},Z)$-rescaling of $\mathcal{A}$ and assume that $\mathcal{A}|_Z/\mathcal{B}$ is a sum of line bundles. Then:
	\begin{itemize}\itemsep=0pt
	\item $w_1([\mathcal{A}\text{:}\mathcal{B}]) = w_1(\mathcal{A}) + k w_1(L_Z)$;
	\item $w_2([\mathcal{A}\text{:}\mathcal{B}]) = w_2(\mathcal{A}) + k w_1(L_Z) \cup w_1(\mathcal{A}) + \frac{k(k-1)}2 w_1(L_Z)^2$;
	\item $p_1([\mathcal{A}\text{:}\mathcal{B}]) = p_1(\mathcal{A})$, if $X$ is orientable and four-dimensional.
	\end{itemize}
\end{Proposition}
\begin{proof} By Proposition~\ref{prop:logbundleiso} we have that $[\mathcal{A}\text{:}\mathcal{B}] \oplus k L_Z \cong \mathcal{A} \oplus k \underline{\mathbb{R}}$ using the shorthand notation $k L = L \oplus \dots \oplus L$ with $k$ copies. This implies using Proposition~\ref{prop:charclasses}.i) by $k$-fold induction that $w([\mathcal{A}\text{:}\mathcal{B}]) \cup (1 + w_1(L_Z))^k = w(\mathcal{A}) \cup 1^k$. We have $(1 + w_1(L_Z))^k \equiv 1 + k w_1(L_Z) + \frac{k(k-1)}2 w_1(L_Z)^2$ up to degree two. In degree one this gives $w_1([\mathcal{A}\text{:}\mathcal{B}]) = w_1(\mathcal{A}) + k w_1(L_Z)$ as desired, while in degree two it instead gives $w_2([\mathcal{A}\text{:}\mathcal{B}]) = w_2(\mathcal{A}) + k w_1(L_Z) \cup w_1(\mathcal{A}) + \frac{k(k-1)}2 w_1(L_Z)^2$. The last property follows because by Proposition~\ref{prop:charclasses}.iii) we have $2p_1([\mathcal{A}\text{:}\mathcal{B}]) = 2p_1(\mathcal{A}) \in H^4(X;\mathbb{Z})$, and the hypothesis ensures that $H^4(X;\mathbb{Z}) \cong \mathbb{Z}$ has no two-torsion (cf.\ Proposition~\ref{prop:logcharclass}).
\end{proof}

\begin{Remark}\label{rem:zeronosymp} We will have no direct use for $w_2(\mathcal{B}_Z)$ and $p_1(\mathcal{B}_Z)$ (when $X$ is four-dimensional), as by \cite[Proposition 2.21]{Lanius16} we know that $\mathcal{B}_Z$ does not admit Lie algebroid symplectic structures when $\dim X \geq 4$. This follows from studying the ring structure of the space $\Omega^\bullet(\mathcal{B}_Z)$.
\end{Remark}

For the $b^k$-tangent bundles we can give an alternate proof to determine the first Stiefel--Whitney class, which does not require the assumption that $k$ is odd.

\begin{Proposition}\label{prop:charclassesbk} Let $(X,Z)$ be a log pair with a defining $(k-1)$-jet $j$ at $Z$. Then we have $w_1\big(\mathcal{A}_Z^k\big) = w_1(TX) + k w_1(L_Z)$.
\end{Proposition}

\begin{proof} From the local description of the Lie algebroid $\mathcal{A}_Z^k$, we see that the anchor map $\rho\colon \mathcal{A}_Z^k \to TX$ leads to a divisor $\big(\det(TX)\otimes\det\big(\mathcal{A}_Z^k\big)^*, \det(\rho)\big)$, which is isomorphic to $(L_Z^k,s^k)$, the $k$th power of a log divisor. Thus $w_1\big(\mathcal{A}_Z^k\big) - w_1(TX) = k w_1(L_Z)$, and hence the conclusion follows.
\end{proof}

We can further determine the first Stiefel--Whitney class of the bundles $\mathcal{A}_{|D|}$ and $\mathcal{A}_W$. This is rather simple, because of the fact that closed submanifolds of codimension two arise.

\begin{Proposition}\label{prop:ellelllogw1} Let $(X,|D|)$ and $(X',W)$ be an elliptic and elliptic-log pair. Then:
	\begin{itemize}\itemsep=0pt
	\item $w_1(\mathcal{A}_D) = w_1(TX)$;
	\item $w_1(\mathcal{A}_W) = w_1(\mathcal{A}_Z) = w_1(TX) + w_1(L_{Z})$, if $I_W = I_{Z} \otimes I_{|D'|}$.
	\end{itemize}
\end{Proposition}

To prove this we first turn to an auxilliary lemma regarding triviality of line bundles.

\begin{Lemma}\label{lem:bundletriviality} Let $L \to X$ be a real line bundle with a section vanishing only on a submanifold of codimension at least two. Then $L$ is trivial, i.e., $w_1(L) = 0 \in H^1(X;\mathbb{Z}_2)$. Consequently, if $(\varphi,{\rm id}_X)\colon E \to F$ is a base-preserving vector bundle morphism which is an isomorphism outside a submanifold of codimension at least two in $X$, then $w_1(E) = w_1(F)$.
\end{Lemma}

\begin{proof} Let $N \subseteq X$ be that submanifold and let $\iota\colon X\backslash N \hookrightarrow N$ be the inclusion. By a standard fact (see \cite[Theorem~2.3, p.~146]{Godbillon71}), the group homomorphism $\iota_*\colon \pi_1(X\backslash N,x) \to \pi_1(X,x)$ is a~surjection between fundamental groups, where $x \in X \backslash N$. This implies by abelianization (using the Hurewicz theorem) and dualizing that the map $\iota^*\colon H^1(X) \to H^1(X\backslash N)$ is an injection. As~$L$ is trivial on $X \backslash N$ by hypothesis, we have $w_1(L|_{X\backslash N}) = 0$, so that $w_1(L) = 0$.

The condition on $\varphi$ being generically an isomorphism implies that ${\rm rank}(E) = {\rm rank}(F)$. Equivalently, using $\det(\varphi)\colon \det(E) \to \det(F)$, the pair $(\det(F) \otimes \det(E)^*, \det(\varphi))$ is a divisor, and $\det(\varphi)$ vanishes only on a submanifold of codimension at least two by hypothesis. The first part then implies that $w_1(\det(F) \otimes \det(E)^*) = 0$, from which the conclusion follows.
\end{proof}

\begin{proof}[Proof of Proposition~\ref{prop:ellelllogw1}] The natural maps $\rho_{\mathcal{A}_{|D|}}\colon \mathcal{A}_{|D|} \to TX$ and $\varphi_{\mathcal{A}_W}\colon \mathcal{A}_W \to \mathcal{A}_Z$ are both isomorphisms outside of $D$ and $D'$ respectively, both of which are of codimension two. Consequently Lemma~\ref{lem:bundletriviality} applies, hence the result follows (using Proposition~\ref{prop:logcharclass} for $\mathcal{A}_Z$).
\end{proof}

\section{Homotopical obstructions}

\subsection{Obstructions from orientability}\label{sec:existanambu}

In this section we discuss orientability for the Lie algebroids $\mathcal{A}_Z^k$, $\mathcal{B}_Z$ and $\mathcal{C}_Z$ associated to log pairs $(X,Z)$, and for the Lie algebroids $\mathcal{A}_{|D|}$ and $\mathcal{A}_W$ given elliptic and elliptic-log pairs $(X,|D|)$ and $(X,W)$. This settles when these Lie algebroids admit symplectic structures in dimension two, and gives an obstruction to their existence in arbitrary dimensions, noting Proposition~\ref{prop:symplaconsequences}. They moreover characterize the existence of $\mathcal{A}$-Nambu structures of highest degree (i.e., nonvanishing sections $\Pi \in \Gamma(\det(\mathcal{A})$).

Given a Lie algebroid $\mathcal{A} \to X$, note that it admitting an orientation does not depend on the Lie algebroid structure of $\mathcal{A}$, and happens if and only if $w_1(\mathcal{A}) = 0$. From our earlier computations we can readily conclude the following (Theorem~\ref{thm:introorientable}):

\begin{Theorem}\label{thm:orientable} Let $\mathcal{A} \to X^n$ be a symplectic Lie algebroid. Then in $H^1(X;\mathbb{Z}_2)$ we have:
	\begin{itemize}\itemsep=0pt
	\item If $\mathcal{A} = \mathcal{A}_Z^k$, then $w_1(TX) + k {\rm PD}_{\mathbb{Z}_2}[Z] = 0$ $($cf.\ {\rm \cite{MirandaPlanas18two,MirandaPlanas18})};
	\item If $\mathcal{A} = \mathcal{B}_Z$, then $w_1(TX) = 0$;
	\item If $\mathcal{A} = \mathcal{C}_Z$, then $w_1(TX) + {\rm PD}_{\mathbb{Z}_2}[Z] = 0$;
	\item If $\mathcal{A} = \mathcal{A}_{|D|}$, then $w_1(TX) = 0$;
	\item If $\mathcal{A} = \mathcal{A}_W$, then $w_1(TX) + {\rm PD}_{\mathbb{Z}_2}[Z] = 0$.
	\end{itemize}
\end{Theorem}

\begin{proof} If the Lie algebroid $\mathcal{A}$ is symplectic, by Proposition~\ref{prop:symplaconsequences} it must be orientable, so that we see that $w_1(\mathcal{A}) = 0$. The result then follows from Proposition~\ref{prop:logcharclass}, Proposition~\ref{prop:charclassesbk} and Proposition~\ref{prop:ellelllogw1}, noting that $n$ is even. We moreover use Lemma~\ref{lem:logdivisor} for the fact that $w_1(L_Z) = {\rm PD}_{\mathbb{Z}_2}[Z] \in H^1(X;\mathbb{Z}_2)$.
\end{proof}

\subsection{Obstructions from complex structures}\label{sec:existenceaacs}

In this section we discuss when some Lie algebroids of interest can admit a complex structure. For this we use Theorem~\ref{thm:wuacs} together with our earlier computations of characteristic classes (see Proposition~\ref{prop:logcharclass}). Due to Proposition~\ref{prop:symplaconsequences} this provides obstructions to when these Lie algebroids can be symplectic.

Let $\big(X^4,Z\big)$ be a four-dimensional log pair with $X$ oriented. Consider a defining $(k-1)$-jet for $Z$ with $k$ odd, and its $b^k$-tangent bundle $\mathcal{A}_Z^k$, which recall includes the log-tangent bundle if $k = 1$. Assume that an orientation for $\mathcal{A}_Z^k$ is given. Then we can define the following:

\begin{Definition}\label{defn:discrepancy} Given orientations on the bundles $\mathcal{A}_Z^k$ and $TX$, the integer $f_k(X,Z)$ of $Z$ is defined as the difference $2 f_k(X,Z) := e\big(\mathcal{A}_Z^k\big) - e(TX) \in H^4(X;\mathbb{Z}) \cong \mathbb{Z}$.
\end{Definition}

\begin{Lemma}\label{lem:kdiscrepancy} In the situation above the integer $f_k(X,Z)$ is well-defined, i.e., the difference in Euler classes of $\mathcal{A}_Z^k$ and $TX$ is even. Further, we have $f_k(X,Z) \equiv f_1(X,Z) \pmod{2}$.
\end{Lemma}

We will henceforth write $f(X,Z) := f_1(X,Z)$.

\begin{proof} Recall that the Euler class of an oriented vector bundle reduces mod $2$ to its top Stiefel--Whitney class. Because both $\mathcal{A}_Z^k$ and $TX$ are oriented, we have $w_1\big(\mathcal{A}_Z^k\big) = w_1(TX) = 0$, hence $w_1(L_Z) = 0$ by Proposition~\ref{prop:logcharclass}. This means that $L_Z$ is trivial, so that by Proposition~\ref{prop:bkbundleiso} we have that $\mathcal{A}_Z^k \oplus \underline{\mathbb{R}} \cong TX \oplus \underline{\mathbb{R}}$. Using Proposition~\ref{prop:charclasses}(i) this implies that $w\big(\mathcal{A}_Z^k\big) = w(TX)$, so that in particular $e\big(\mathcal{A}_Z^k\big) \equiv e(TX) \pmod{2}$ as desired.

For the second statement, we remark that there is a more geometric description of $f_k(X,Z)$. If $\mathcal{A}^k_Z$ is oriented and $X$ is orientable, any choice of orientation for $TX$ does not agree with the orientation on the isomorphism locus $X \backslash Z$ induced by $\mathcal{A}_Z$ if and only if $k$ is odd. This follows from the local description of the bundle $\mathcal{A}_Z^k$ (cf.\ \cite{CavalcantiKlaasse19} for when $k = 1$). If $k$ is odd, because $Z$ is then separating due to Proposition~\ref{prop:logcharclass}, after a choice of orientation on $TX$, we can write $X \backslash Z = X_+ \sqcup X_-$, where $X_\pm$ denote the subsets where the orientations from $TX$ and $\mathcal{A}_Z^k$ do or do not agree. Then we see that
\begin{displaymath}
	f_k(X,Z) = - \langle e(TX), [X_-] \rangle = - \chi(X_-).
\end{displaymath}
We see here that the right-hand side does not depend on $k$, nor does the decomposition of $X\backslash Z$ into~$X_\pm$. It follows that in fact $f_k(X,Z) = f_1(X,Z)$ if $k$ is odd.
\end{proof}

We can now state our obstruction to the existence of an $\mathcal{A}^k_Z$-almost-complex structure.

\begin{Theorem}\label{thm:azalmostcplx} Let $\big(X^4,Z\big)$ be a compact oriented $\mathcal{A}^k_Z$-almost-complex log pair, with $k$ odd. Then we have $\big\langle \big[c_1^2\big(\mathcal{A}^k_Z\big)\big], [X] \big\rangle = 3 \sigma(X) + 2 \chi(X) + 4 f_k(X,Z)$, and $b_2^+(X) + b_1(X) + f_k(X,Z)$ is odd.
\end{Theorem}

Here $\chi(X)$ is the Euler characteristic, and $\sigma(X) = b_2^+(X) - b_2^-(X)$ is the signature of~$X$.
The following proof is similar to the case when $Z = \varnothing$, see \cite[Theorem~1.4.13]{GompfStipsicz99}.

\begin{proof} By Theorem~\ref{thm:wuacs} we have, using the definition of $f_k(X,Z)$, that
\begin{displaymath}
	c_1^2\big(\mathcal{A}_Z^k\big) = p_1\big(\mathcal{A}_Z^k\big) + 2 e\big(\mathcal{A}_Z^k\big) = p_1\big(\mathcal{A}_Z^k\big) + 2 e(TX) + 4 f_k(X,Z).
\end{displaymath}
Using the other part of Theorem~\ref{thm:wuacs} together with Proposition~\ref{prop:logcharclass} we get (as $w_1(L_Z) = 0$)
\begin{displaymath}
	c_1\big(\mathcal{A}_Z^k\big) \!\!\!\!\mod 2 \equiv w_2\big(\mathcal{A}_Z^k\big) = w_2(TX) \in H^2(X;\mathbb{Z}_2),
\end{displaymath}
so that $c_1\big(\mathcal{A}_Z^k\big)$ is \emph{characteristic}, i.e., it reduces modulo two to $w_2(TX)$. By Van der Blij's lemma \cite[Lemma~1.2.20]{GompfStipsicz99} we obtain $c_1^2\big(\mathcal{A}_Z^k\big) \equiv \sigma(X) \!\!\mod 8$. Using Proposition~\ref{prop:logcharclass} again we have $p_1\big(\mathcal{A}_Z^k\big) = p_1(TX)$, which integrates to $3 \sigma(X)$ by the Hirzebruch signature theorem. This implies that $\sigma(X) + \chi(X) + 2 f_k(X,Z) \equiv 0 \!\!\mod 4$, from which the conclusion follows.
\end{proof}

Out of this we can draw the following consequence (the first part of Theorem~\ref{thm:intrologsymp}).

\begin{Theorem}\label{thm:azsympstr} Let $\big(X^4,Z\big)$ be a compact oriented $b^k$-symplectic manifold, with $k$ odd. Then:
	\begin{displaymath}
		b_2^+(X) + b_1(X) + f(X,Z) \quad \text{is odd}.
	\end{displaymath}
\end{Theorem}

\begin{proof} Because $\mathcal{A}_Z^k$ is a symplectic Lie algebroid, it is also complex by Proposition~\ref{prop:symplaconsequences}. Due to this $(X,Z)$ has an $\mathcal{A}_Z^k$-almost-complex structure, so that the result follows from Theorem~\ref{thm:azalmostcplx} and applying Lemma~\ref{lem:kdiscrepancy} to replace $f_k(X,Z)$ by $f(X,Z)$ after reducing modulo two.
\end{proof}

We would like to stress again that this obstruction, if $k$ is even, agrees with that for $X$ to admit an almost-complex structure in the usual sense (see \cite[Theorem~1.4.13]{GompfStipsicz99}).
A similar thing can be done for the scattering tangent bundle $\mathcal{C}_Z$, as we now discuss. Consider again a four-dimensional oriented log pair $\big(X^4, Z\big)$, and assume that an orientation for $\mathcal{C}_Z$ is given. There is an analogous difference in Euler classes here, similar to Definition~\ref{defn:discrepancy}.

\begin{Definition} Given orientations on the bundles $\mathcal{C}_Z$ and $TX$, the integer $f_{\rm sc}(X,Z)$ of $Z$ is defined as the difference $2 f_{\rm sc}(X,Z) := e(\mathcal{C}_Z) - e(TX) \in H^4(X;\mathbb{Z}) \cong \mathbb{Z}$.
\end{Definition}

In fact, we can quickly relate the integer $f_{\rm sc}(X,Z)$ to $f(X,Z)$ defined previously.

\begin{Lemma}\label{lem:scdiscrepancy} Let $\big(X^{2n}, Z\big)$ be a log pair, and choose orientations on $\mathcal{C}_Z$ and $TX$. Then $\mathcal{A}_Z$ is naturally oriented, and we have the equality $f_{\rm sc}(X,Z) \equiv f(X,Z) \pmod{2}$.
\end{Lemma}

\begin{proof} The natural Lie algebroid morphism $\varphi\colon \mathcal{C}_Z \to \mathcal{A}_Z$ can be used to orient $\mathcal{A}_Z$. Note that because the dimension of $X$ is even, we have that $w_1(\mathcal{C}_Z) = w_1(\mathcal{A}_Z)$ by Proposition~\ref{prop:logcharclass}. From the fact that the divisor ideal of $\varphi$ is given by $I_\varphi = I_Z^{2n}$, or by the local description of $\mathcal{C}_Z$, it readily follows that as for $\mathcal{A}_Z$, the orientation on $X\backslash Z$ induced by $\mathcal{C}_Z$ similarly cannot match the one induced from $TX$ everywhere, from which the result follows.
\end{proof}

Using this we can obtain an obstruction to the existence of a $\mathcal{C}_Z$-almost-complex structure.

\begin{Theorem}\label{thm:scatteringalmostcplx} Let $\big(X^4,Z\big)$ be a compact oriented $\mathcal{C}_Z$-almost-complex log pair. Then we have $\big\langle \big[c_1^2(\mathcal{C}_Z)\big], [X] \big\rangle = 3 \sigma(X) + 2 \chi(X) + 4 f_{\rm sc}(X,Z)$, and $b_2^+(X) + b_1(X) + f(X,Z)$ is odd.
\end{Theorem}

\begin{proof} The proof follows the along the same lines as that of Theorem~\ref{thm:azalmostcplx}. The first statement follows by definition of $f_{\rm sc}(X,Z)$. Further, here Proposition~\ref{prop:logcharclass} provides
\begin{displaymath}
	c_1(\mathcal{C}_Z) \!\!\mod 2 \equiv w_2(\mathcal{C}_Z) = w_2(TX),
\end{displaymath}
showing again that $c_1(\mathcal{C}_Z)$ is characteristic. Proposition~\ref{prop:logcharclass} also gives $p_1(\mathcal{C}_Z) = p_1(TX)$, which implies that $\sigma(X) + \chi(X) + 2 f_{\rm sc}(X) \equiv 0 \!\!\mod 4$. By Lemma~\ref{lem:scdiscrepancy} we can replace~$f_{\rm sc}(X,Z)$ by~$f(X,Z)$ after reduction modulo two.
\end{proof}

As for Theorem~\ref{thm:azsympstr}, we can draw the following consequence (second part of Theorem~\ref{thm:intrologsymp}).

\begin{Theorem}\label{thm:scatteringsymp} Let $\big(X^4,Z\big)$ be a compact oriented scattering-symplectic manifold. Then:
	\begin{displaymath}
		b_2^+(X) + b_1(X) + f(X,Z) \quad \text{is odd}.
	\end{displaymath}
\end{Theorem}

In other words, we see that the obstruction for scattering symplectic structures obtained via the existence of almost-complex structures is identical to that of log-symplectic structures.

\begin{proof} Because $\mathcal{C}_Z$ is a symplectic Lie algebroid, it is also complex by Proposition~\ref{prop:symplaconsequences}. Due to this $(X,Z)$ has a $\mathcal{C}_Z$-almost-complex structure, hence the result follows from Theorem~\ref{thm:scatteringalmostcplx}.
\end{proof}

\begin{Remark} One wonders whether Proposition~\ref{prop:symplaconsequences} can be used effectively for other symplectic Lie algebroids in dimension four, for example the elliptic tangent bundle $\mathcal{A}_{|D|}$. Note that elliptic symplectic structures (of zero elliptic residue) can exist on $\mathcal{A}_{|D|}$ both in cases when $X$ is and is not almost-complex (cf.\ \cite{BehrensCavalcantiKlaasse18}), depending on the coorientability of $D$ as measured by $w_1(ND) \in H^1(D;\mathbb{Z}_2)$. We see there is nontrivial dependence on the locus $D$ in this case.
\end{Remark}

To illustrate Theorem~\ref{thm:azalmostcplx}, we determine the parity of $f(X,Z)$ in the following situation. As is explained in the proof of Lemma~\ref{lem:kdiscrepancy}, if both $X$ and $\mathcal{A}_Z$ are oriented, then $Z$ must be separating and decompose $X\backslash Z = X_+ \sqcup X_-$ according to whether the orientations agree.

\begin{Corollary}\label{cor:aznosympconn} Let $(X,Z)$ be a compact oriented four-dimensional log pair which is $\mathcal{A}_Z$-almost-complex, such that $X$ is not almost-complex. Then $f(X,Z)$ is odd, and the log pair $\big(X_- \cup \big(X_+ \# \mathbb{C} P^2\big), Z\big)$ does not admit an $\mathcal{A}_Z$-symplectic structure.
\end{Corollary}

\begin{proof} If $X$ is not almost-complex, then $b_2^+(X) + b_1(X) \equiv 0 \pmod 2$, while because $(X,Z)$ is $\mathcal{A}_Z$-almost-complex we obtain from Theorem~\ref{thm:azalmostcplx} that $b_2^+(X) + b_1(X) + f(X,Z) \equiv 1 \pmod 2$. We conclude that $f(X,Z) \equiv 1 \pmod 2$. If we perform a connected sum with $\mathbb{C} P^2$ in the subset~$X_+$ to form the manifold $X' = X_- \cup \big(X_+ \# \mathbb{C} P^2\big)$, we see that $b_2^+(X') = b_2^+(X) + 1$ while $f(X',Z) = f(X,Z)$. Hence then $b_2^+(X') + b_1(X') + f(X',Z) \equiv 0 \pmod 2$, so that by applying Theorem~\ref{thm:azsympstr} we see that $(X',Z)$ does not admit an $\mathcal{A}_Z$-symplectic structure.
\end{proof}

We finish by giving a simple example of how to apply the above results.

\begin{Example}[$3 \mathbb{C} P^2 \# \overline{\mathbb{C} P}^2$] The manifold $X = 2 \mathbb{C} P^2 \# \overline{\mathbb{C} P}^2$ admits a log-symplectic structure with $Z = S^1 \times S^2$ (see \cite{Cavalcanti17}), and $b_2^+(X) = 2$ and $b_1(X) = 0$. Hence $X$ is not almost-complex while $(X,Z)$ is $\mathcal{A}_Z$-almost-complex, and $f(X,Z)$ is odd. By Corollary~\ref{cor:aznosympconn} the manifold $X' = 3 \mathbb{C} P^2 \# \overline{\mathbb{C} P}^2$ does not admit a log-symplectic structure with the given $Z = S^1 \times S^2$. Note, however, that by~\cite{Cavalcanti17} again, it does admit a log-symplectic structure with degeneracy locus diffeomorphic to $S^1 \times S^2$, but this must necessarily be a different hypersurface.
\end{Example}

\begin{Remark} It seems somewhat nontrivial to determine the integer $f(X,Z)$ of a separating log pair, or even just its parity. We are only able to do so indirectly, cf.\ Corollary~\ref{cor:aznosympconn}.
\end{Remark}

\subsection*{Acknowledgements}

This work is partially based on \cite[Chapter 11]{Klaasse17}, and was supported by ERC consolidator grant 646649 ``SymplecticEinstein'', and by VIDI grant 639.032.221 from NWO, the Netherlands Organisation for Scientific Research. The author would like to thank Gil Cavalcanti for useful discussions, and the anonymous referees for their suggestions.

\pdfbookmark[1]{References}{ref}
\LastPageEnding

\end{document}